\documentclass[12pt]{amsart}

\newtheorem{theorem}{Theorem}[section]
\newtheorem{lemma}[theorem]{Lemma}
\newtheorem{proposition}[theorem]{Proposition}
\newtheorem{corollary}[theorem]{Corollary}
\newtheorem{problem}[theorem]{Problem}
\newtheorem{conjecture}[theorem]{Conjecture}

\theoremstyle{definition}

\newtheorem{example}[theorem]{Example}
\newtheorem{remark}[theorem]{Remark}

\theoremstyle{remark}

\usepackage{amscd,amssymb}

\begin{document}
\baselineskip 15pt

\title[The Strong Anick Conjecture Is True]
{The strong Anick conjecture is true}

\author[Vesselin Drensky and Jie-Tai Yu]
{Vesselin Drensky and Jie-Tai Yu}
\address{Institute of Mathematics and Informatics,
Bulgarian Academy of Sciences, 1113 Sofia, Bulgaria}
\email{drensky@math.bas.bg}
\address{Department of Mathematics, The University of Hong Kong, Hong
Kong SAR, CHINA}
\email{yujt@hkucc.hku.hk}

\thanks
{The research of Vesselin Drensky was partially supported by Grant
MM-1106/2001 of the Bulgarian National Science Fund.}
\thanks{The research of Jie-Tai Yu was partially
supported by a Hong Kong RGC-CERG Grant.}

\subjclass[2000]{Primary 16S10. Secondary 16W20; 13B10; 13B25.}
\keywords{Automorphisms of free and polynomial algebras, tame
automorphisms,
wild automorphisms, coordinates in free algebras}

\begin{abstract}
Recently Umirbaev has proved the long-standing Anick conjecture,
that is, there exist wild automorphisms of the free associative
algebra $K\langle x,y,z\rangle$ over a field $K$ of characteristic
0. In particular, the well-known Anick automorphism is wild. In
this article we obtain a stronger result (the Strong Anick
Conjecture that implies the Anick Conjecture). Namely, we prove
that there exist wild coordinates of $K\langle x,y,z\rangle$. In
particular, the two nontrivial coordinates in the Anick
automorphism are both wild. We establish a similar result for
several large classes of automorphisms of $K\langle x,y,z\rangle$.
We also find a large new class of wild automorphisms of
$K\langle x,y,z\rangle$ which is not covered by the results
of Umirbaev. Finally, we study
the lifting problem for automorphisms and coordinates of
polynomial algebras, free metabelian algebras and free associative
algebras and obtain some interesting new results.
\end{abstract}
\maketitle

\section{Introduction}

Let $K$ be a field of characteristic 0 and let
$X=\{x_1,\ldots,x_n\}$, $n\geq 2$, be a finite set. We denote by
$K[X]$ the polynomial algebra in the set of variables $X$ and by
$K\langle X\rangle$ the free associative algebra (or the algebra
of polynomials in the set $X$ of noncommuting variables)
with the same set of free generators.
For small $n$ we shall denote the free generators also with $x,y$, etc.
We write the automorphisms of $K[X]$ and $K\langle X\rangle$
as $n$-tuples of the images of the coordinates,
and $\varphi=(f_1,\ldots,f_n)$
means that $\varphi(x_j)=f_j(X)=f_j(x_1,\ldots,x_n)$,
$j=1,\ldots,n$.
Also, the product $\varphi\psi$ of the automorphisms $\varphi$ and
$\psi$
acts on $u(X)$ by $(\varphi\psi)(u)=\varphi(\psi(u))$.
An automorphism of $K[X]$ or $K\langle X\rangle$
is called {\it elementary}, if it is of the form
\[
(x_1,\ldots,x_{j-1},\alpha
x_j+f(x_1,\ldots,x_{j-1},x_{j+1},\ldots,x_n),
x_{j+1},\ldots,x_n),\quad\alpha\in K^{\ast},
\]
and the polynomial $f(x_1,\ldots,x_{j-1},x_{j+1},\ldots,x_n)$ does not
depend on the
variable $x_j$. The automorphisms belonging to the group generated
by elementary automorphisms are called {\it tame}, otherwise they are
{\it wild}. A polynomial $p\in K[X]$ is called a {\it coordinate} if it
is an automorphic image of $x_1$. Moreover, a coordinate $p\in K[X]$ is
called
{\it tame} if there exists a tame automorphism
$\varphi\in\text{\rm Aut}K[X]$ such that $\varphi(x_1)=p$, otherwise
$p$ is called a {\it wild coordinate}. One defines in a similar way the
coordinates, tame and wild coordinates
of $K\langle X\rangle$ and other relatively free algebras. In
noncommutative
algebra coordinates are often called $primitive\ elements$.
Obviously, the existence of wild coordinates  implies the existence of
wild automorphisms, but
not vice versa in general.

Problems concerning automorphisms of free objects are
similar for free groups, polynomial algebras, free associative and free
Lie
algebras, for relatively free groups and algebras,
see the recent book \cite{MSY} by Mikhalev, Shpilrain, and Yu.
One of the most important problems in the theory of automorphisms
of polynomial and free associative algebras is
whether all automorphisms of $K[X]$ and $K\langle X\rangle$ are tame.
The answer is in
the affirmative for $K[x,y]$ (Jung--van der Kulk \cite{J, K}) and for
$K\langle x,y\rangle$
(Czerniakiewicz--Makar-Limanov \cite{Cz, ML1, ML2}). In 1970,
Nagata \cite{N} constructed his famous
automorphism of $K[x,y,z]$ which is wild as an automorphism fixing $z$
and conjectured that it is wild also as a usual automorphism of
$K[x,y,z]$.
More than 30 years later Shestakov and Umirbaev \cite{SU1, SU2, SU3},
using methods of Poisson brackets, degree estimations and peak
reduction,
proved the Nagata Conjecture, in particular, they proved that the
Nagata automorphism is wild.
Umirbaev and J.-T. Yu \cite{UY} proved the Strong Nagata Conjecture,
namely, there exist wild coordinates of
$K[x,y,z]$. In particular, the two nontrivial Nagata coordinates are
both wild.
There were also attempts, unfortunately unsuccessful,
to lift the Nagata automorphism to an automorphism of
$K\langle x,y,z\rangle$, see, for instance, our paper with Gutierrez
\cite{DGY}.
(A lifting of the Nagata automorphism would provide immediately a wild
automorphism
of $K\langle x,y,z\rangle$.) There is another automorphism of $K\langle
x,y,z\rangle$,
suggested by Anick,
\[
(x+y(xy-yz),y,z+(xy-yz)y)\in \text{Aut}K\langle x,y,z\rangle,
\]
see the book by Cohn \cite{C2}, p.~343, which
was suspected to be wild, although it induces a tame automorphism of
$K[x,y,z]$. Exchanging the places of $y$ and $z$ in the above
automorphism, we replace it
with the automorphism
\[
\omega=(x+z(xz-zy),y+(xz-zy)z,z).
\]
It has the property that fixes $z$ and $\omega(x)$ and $\omega(y)$ are
linear
in $x$ and $y$. From now on, we shall refer to $\omega$ as the {\it
Anick automorphism}.

\begin{conjecture} {\bf (Anick Conjecture)}
There exist wild automorphisms in $\text{\rm Aut}K\langle
x,y,z\rangle$.
In particular, the Anick automorphism is wild.
\end{conjecture}

We established \cite{DY3} that
the Anick automorphism is wild in the class of automorphisms fixing the
variable $z$,
and very recently Umirbaev \cite{U2} has proved the Anick Conjecture,
that is,
the wildness of the Anick automorphism as an automorphism of $K\langle
x,y,z\rangle$.

The work of Nagata \cite{N} motivated the study of automorphisms fixing
a variable.
More generally, we introduce another set $Z=\{z_1,\ldots,z_m\}$ and
consider the algebras
$K[X,Z]$ and $K\langle X,Z\rangle$, freely generated by the set $X\cup
Z$.
Studying automorphisms fixing the set $Z$, we use the same notation
$\varphi=(f_1,\ldots,f_n,z_1,\ldots,z_m)$, meaning that
$\varphi(x_j)=f_j(X,Z)$,
$j=1,\ldots,n$, and $\varphi(z_i)=z_i$, $i=1,\ldots,m$,
and denote the corresponding automorphism group by
$\text{\rm Aut}_ZK[X,Z]$ and $\text{\rm Aut}_ZK\langle X,Z\rangle$,
respectively.
Again, an automorphism is $Z$-elementary if it is
of the form $\varphi=(x_1,\ldots,\alpha
x_j+f(X,Z),\ldots,x_n,z_1,\ldots,z_m)$,
and $f(X,Z)$ does not depend on $x_j$. The automorphism is tame in the
class of
automorphisms fixing $Z$ (or $Z$-tame) if it belongs to the group
generated by $Z$-elementary
automorphisms. In the general case, this group may be smaller than the
group generated by
$Z$-triangular and $Z$-affine automorphisms because some $Z$-affine
automorphisms may be
not products of $Z$-elementary ones.

An important consequence of  \cite{SU1,SU2, SU3, DY1, DY2} is the
result
that every $z$-wild automorphism of $K[x,y,z]$ is a wild automorphism
of $K[x,y,z]$.
A way to construct a large class of
such automorphisms was given by us in \cite{DY2}.

The next step of studying automorphisms of free algebras
is to study {\it coordinates}, or automorphic images of $x_1$. In
noncommutative algebra
coordinates are also called {\it primitive elements}. There are
algorithms to recognize
coordinates of $K[x,y]$ (Shpilrain and Yu \cite{SY}) and
$z$-coordinates of $K[x,y,z]$, see
our paper \cite{DY2} as well as our survey \cite{DY1}.
As a consequence of their proof of the Strong Nagata Conjecture,
Umirbaev and J.-T. Yu \cite {UY} established that if
$f(x,y,z)$ is a $z$-wild coordinate in $K[x,y,z]$, then it is
immediately a wild coordinate in $K[x,y,z]$.

As a common sense, results for (commutative) polynomial algebras
inspire problems on free
associative algebras as there is a natural surjective homomorphism
from $K\langle X\rangle$ to $K[X]$.
which induces a natural (not necessarily surjective) homomorphism
from $\text {Aut}K\langle X\rangle$ to $\text{Aut}K[X]$.
In this paper we shall be interested in the following problem
motivated by \cite{UY}:

\begin{problem}\label{wild coordinates}
If $f(X,Z)\in K\langle X,Z\rangle$ is an image of $x_1$ under a
$Z$-wild
automorphism, is there a tame automorphism (maybe not fixing $Z$)
which also sends $x_1$ to $f(X,Z)$?
\end{problem}

If such a tame automorphism does not exist, then we call $f(X,Z)$
a {\it wild coordinate} of $K\langle X,Z\rangle$.

As an analog of the Strong Nagata Conjecture in \cite{UY}, we state

\begin{conjecture} {\bf (Strong Anick Conjecture)}
There exist wild coordinates in $K\langle x,y,z\rangle$.
In particular, the two nontrivial coordinates of the Anick automorphism
are both wild.
\end{conjecture}
The Anick automorphism has the property that it fixes $z$ and
$\omega(x)$
and $\omega(y)$ are linear
in $x$ and $y$. In our paper \cite{DY3} we showed that such an
automorphism
of $K\langle x,y,z\rangle$ is $z$-tame if and only if certain $2\times
2$ matrix
with entries from $K[z_1,z_2]$ is a product of elementary matrices.
This idea has been used by Umirbaev \cite{U2} in the final step of his
proof of the
the wildness of the Anick automorphism.
Now we show that if $f(x,y,z)$ is a wild $z$-coordinate in $K\langle
x,y,z\rangle$,
and $f(x,y,z)$ is linear in $x,y$, then $f(x,y,z)$ is also wild in the
sense
of Problem \ref{wild coordinates}. This is one of the main results of
the paper.
It immediately gives an affirmative
answer to the Strong Anick Conjecture.

The class of wild automorphisms of $K\langle x,y,z\rangle$
discovered by Umirbaev \cite{U2}
is larger than the class of $z$-wild automorphisms $(f,g,z)$ such that
the polynomials $f,g$ are linear in $x,y$. Our method also gives that
all automorphisms of the class of Umirbaev have the property that at
least two
of their coordinates are wild. The same result holds for another large
class of
automorphisms of $K\langle x,y,z\rangle$ which is not covered by
Umirbaev \cite{U2}.

Our main result suggests an algorithm  deciding whether a polynomial
$f(x,y,z)\in K\langle x,y,z\rangle$
which is linear in $x$ and $y$, is a tame coordinate. If it is, then
the algorithm shows how
to find a product of $z$-elementary automorphisms which sends $x$ to
$f(x,y,z)$.
(Of course, in all algorithmic considerations we assume that the ground
field $K$ is
constructive, and we may perform calculations there.)
In this part of the paper we use the approach and the results of
Umirbaev
\cite{U2}, combined with our approach from \cite{DY3}.

On the other hand, we show that the situation is completely different
in the case of
the free metabelian algebra $M(x,y,z)$. We construct an automorphism
which
fixes $y$ and $z$ and cannot be lifted to an automorphism of $K\langle
x,y,z\rangle$.
The proof is based on a test recognizing  some classes of endomorphisms
which are not automorphisms. This test is originally  from group
theory, see Bryant, Gupta, Levin and
Mochizuki \cite{BGLM} and was addapted to algebras by Bryant and
Drensky \cite{BD}.

In addition, we show that an automorphism of $K[X,Z]$
or $K\langle X,Z\rangle$ which is $Z$-wild cannot be lifted to a
$Z$-automorphism of
the absolutely free algebra $K\{X,Z\}$. (As a consequence of a result
of Kurosh \cite{Ku},
all automorphisms of the absolutely free algebra are tame.)
This is equivalent to the fact that there exist no
$Z$-wild automorphisms of $K\{X,Z\}$.

\section{Proof of main results}

Dicks and Lewin \cite{DL} introduced the Jacobian matrix of an
endomorphism
of $K\langle X\rangle$. This is an $n\times n$ matrix with entries from
the tensor
product $K\langle X\rangle\otimes_KK\langle X\rangle^{\text{op}}$ of
the free algebra $K\langle X\rangle$ and its opposite algebra (or its
anti-isomorphic algebra)
$K\langle X\rangle^{\text{op}}$. For $n=2$ they proved that the
Jacobian matrix is invertible over
$K\langle x,y\rangle\otimes_KK\langle x,y\rangle^{\text{op}}$
if and only if the endomorphism is an automorphism. The general case of
any $n$
was established by Schofield \cite{Sc}, which is the Jacobian
Conjecture for free associative algebras.
The partial derivatives and the Jacobian matrix
of Dicks and Lewin can be defined as follows:
\[
\frac{\partial x_i}{\partial x_i}=1,
\quad \frac{\partial x_j}{\partial x_i}=0,\,j\not=i,
\]
and, for  a monomial
$w=x_{i_1}\cdots x_{i_m}\in K\langle X\rangle$
\[
\frac{\partial w}{\partial x_i}=
\sum_{k=1}^m(x_{i_1}\cdots x_{i_{k-1}})\otimes
(x_{i_{k+1}}\cdots x_{i_m})
\frac{\partial x_{i_k}}{\partial x_i},
\]
where $x_{i_1}\cdots x_{i_{k-1}}\in K\langle X\rangle$ and
$x_{i_{k+1}}\cdots x_{i_m}\in K\langle X\rangle^{\text{op}}$.
Then, as usual,
\[
J(\varphi)=\left(\frac{\partial\varphi(x_j)}{\partial x_i}\right),\quad
\varphi\in\text{End}K\langle X\rangle.
\]
We need the following lemma.

\begin{lemma}\label{one nontrivial coordinate}
The only automorphisms of $K\langle X\rangle$  fixing $x_2,\ldots,x_n$
are the tame automorphisms of the form
\[
\tau=(\alpha x_1+f(x_2,\ldots,x_n),x_2,\ldots,x_n),\quad \alpha\in
K^{\ast},
\,f(x_2,\ldots,x_n)\in K\langle x_2,\ldots,x_n\rangle.
\]
\end{lemma}

\begin{proof}
The shortest way to establish the lemma is to use the invertibility
of the Jacobian matrix.
Let $\tau=(g(X),x_2,\ldots,x_n)\in\text{Aut}K\langle X\rangle$
fix $x_2,\ldots,x_n$. Then the matrix
\[
J(\tau)=\left(\begin{matrix}
\frac{\partial g}{\partial x_1}&0&\ldots&0\\
\frac{\partial g}{\partial x_2}&1&\ldots&0\\
\vdots&\vdots&\ddots&\vdots\\
\frac{\partial g}{\partial x_n}&0&\ldots&1\\
\end{matrix}\right)
\]
is invertible over $K\langle X\rangle\otimes_KK\langle
X\rangle^{\text{op}}$
and this implies that $\partial g/\partial x_1$ is equal to a nonzero
constant
$\alpha$. Hence the only term of $g(X)$ depending on $x_1$ is $\alpha
x_1$.
\end{proof}

For $K\langle x,y,z\rangle$, the endomorphisms which fix $z$ and are
linear in
$x$ and $y$ are of the form $\rho=(f(x,y,z),g(x,y,z),z)$, where
\[
f(x,y,z)=\sum_{p,q\geq 0}\alpha_{pq}z^pxz^q+\sum_{p,q\geq
0}\beta_{pq}z^pyz^q+f_0(z),
\]
\[
g(x,y,z)=\sum_{p,q\geq 0}\gamma_{pq}z^pxz^q+\sum_{p,q\geq
0}\delta_{pq}z^pyz^q+g_0(z),
\]
$\alpha_{pq},\beta_{pq},\gamma_{pq},\delta_{pq}\in K$,
and $f_0(z),g_0(z)$ are polynomials in $z$.
Applying the Jacobian matrix of Dicks and Lewin in this concrete case,
in
\cite{DY3} we obtained:

\begin{proposition}\label{linear tame automorphisms fixing variable}
{\rm (i)} The endomorphism $\rho=(f(x,y,z),g(x,y,z),z)$ which fixes $z$
and is linear in $x$ and $y$ is an automorphism if and only if
the $2\times 2$ matrix
\[
J_z(\rho)=\left(\begin{matrix}
\sum_{p,q\geq 0}\alpha_{pq}z_1^pz_2^q&\sum_{p,q\geq
0}\gamma_{pq}z_1^pz_2^q\\
\sum_{p,q\geq 0}\beta_{pq}z_1^pz_2^q&\sum_{p,q\geq
0}\delta_{pq}z_1^pz_2^q\\
\end{matrix}\right)
\]
with entries from $K[z_1,z_2]$ is invertible. All such automorphisms
induce $z$-tame
automorphisms of $K[x,y,z]$.

{\rm (ii)} The automorphism $\rho$ is $z$-tame
if and only if the matrix $J_z(\rho)$
belongs to the group generated by elementary matrices with entries
from $K[z_1,z_2]$.
\end{proposition}

For example, for the Anick automorphism,
\[
J_z(\omega)=\left(\begin{matrix}
1+z_1z_2&z_2^2\\
-z_1^2&1-z_1z_2\\
\end{matrix}\right)
\]
and by a result of Cohn \cite{C1}, the matrix
$J_z(\omega)$ cannot be presented as a product of elementary matrices.

Let $Z=\{z_1,\ldots,z_m\}$. We denote by $GE_2(K[Z])$ the subgroup of
$GL_2(K[Z])$ generated by the diagonal and by the elementary matrices
\[
\left(\begin{matrix}
\alpha_1&0\\
0&\alpha_2\\
\end{matrix}\right),\quad
\left(\begin{matrix}
1&f(Z)\\
0&1\\
\end{matrix}\right),\quad
\left(\begin{matrix}
1&0\\
f(Z)&1\\
\end{matrix}\right)
\]
with entries from $K[Z]$.
There is an algorithm  deciding whether a matrix
in $GL_2(K[Z])$ belongs to $GE_2(K[Z])$.
It was suggested by Tolhuizen, Hollmann, and Kalker \cite{THK}
for the partial ordering by degree and, independently, by Park
\cite{P1, P2}
for any monomial ordering on $K[Z]$. One applies Gaussian
elimination process on the matrix based on the Euclidean division
algorithm
for $K[Z]$. The matrix belongs to $GE_2(K[Z])$
if and only if this procedure brings it to the identity matrix.
For our purposes, we need the following version of the Euclidean
algorithm.
If $a(Z),b(Z)$ are two nonzero polynomials with homogeneous components
of maximal
degree $\overline{a(Z)},\overline{b(Z)}$, respectively, then the
Euclidean algorithm can be applied to $a(Z)$ and $b(Z)$
if $\overline{a(Z)}=\overline{b(Z)}q(Z)$ for some $q(Z)\in K[Z]$
(or $\overline{b(Z)}=\overline{a(Z)}q(Z)$) when we replace
$a(Z)$ with $a(Z)-b(Z)q(Z)$ (or, respectively, we replace $b(Z)$ with
$b(Z)-a(Z)q(Z)$).
In matrix form, these operations correspond, respectively, to
\begin{equation}\label{matrix form A of Euclidean algorithm}
\left(\begin{matrix}
a(Z)-b(Z)q(Z)\\
b(Z)\\
\end{matrix}\right)
=\left(\begin{matrix}
1&-q(Z)\\
0&1\\
\end{matrix}\right)
\left(\begin{matrix}
a(Z)\\
b(Z)\\
\end{matrix}\right),
\end{equation}
\begin{equation}\label{matrix form B of Euclidean algorithm}
\left(\begin{matrix}
a(Z)\\
b(Z)-a(Z)q(Z)\\
\end{matrix}\right)
=\left(\begin{matrix}
1&0\\
-q(Z)&1\\
\end{matrix}\right)
\left(\begin{matrix}
a(Z)\\
b(Z)\\
\end{matrix}\right).
\end{equation}
For us, the most convenient form of the result in
\cite{P1, P2, THK} is as stated in \cite{THK}.

\begin{proposition}\label{Euclidean algorithm for elementary matrices}
Let $a(Z),b(Z)$ be two polynomials in $K[Z]$. Then there exist
$c(Z),d(Z)\in K[Z]$
such that the matrix
\[
G=\left(\begin{matrix}
a(Z)&c(Z)\\
b(Z)&d(Z)\\
\end{matrix}\right)
\]
belongs to $GE_2(K[Z])$ if and only if we can bring the pair
$(a(Z),b(Z))$ to $(\alpha,0)$, $0\not=\alpha\in K$, using the Euclidean
algorithm only.
\end{proposition}

Clearly, in this case the equations (\ref{matrix form A of Euclidean
algorithm})
and (\ref{matrix form B of Euclidean algorithm})
give the decomposition of $G$ as a product of elementary matrices.

We need a description of the free metabelian associative algebra
and a short exposition of the results of Umirbaev.
Recall that the free metabelian algebra
\[
M(X)=K\langle X\rangle/([t_1,t_2][t_3,t_4])^T
\]
is the relatively free algebra of rank $n$ in the variety of
associative algebras
defined by the polynomial identity $[t_1,t_2][t_3,t_4]=0$. In order to
define
partial derivatives and the Jacobian matrix of an endomorphism, we need
two more
sets of variables $U=\{u_1,\ldots,u_n\}$ and $V=\{v_1,\ldots,v_n\}$ of
the same cardinality
as $X$. We consider the polynomial algebra
$K[U,V]$. Changing a little the notation of Umirbaev \cite{U1},
we define formal partial derivatives
$\partial_M/\partial_Mx_i$ assuming that
\[
\frac{\partial_Mx_i}{\partial_Mx_i}=1,
\quad \frac{\partial_Mx_j}{\partial_Mx_i}=0,\,j\not=i,
\]
and, for  a monomial
$w=x_{i_1}\cdots x_{i_m}\in M(X)$
\[
\frac{\partial_Mw}{\partial_Mx_i}=
\sum_{k=1}^mu_{i_1}\cdots u_{i_{k-1}}
v_{i_{k+1}}\ldots v_{i_m}
\frac{\partial_Mx_{i_k}}{\partial_Mx_i}.
\]
These are the homomorphic images of the partial derivatives of
Dicks and Lewin under the natural homomorphism
$K\langle X\rangle\otimes_KK\langle X\rangle^{\text{op}}\to K[U,V]$
which sends $x_i\otimes 1$ and $1\otimes x_j$ to $u_i$ and $v_j$,
respectively.
A polynomial $f(X)\in M(X)$ belongs to the commutator ideal of $M(X)$,
i.e., to the kernel of the natural homomorphism
$M(X)\to K[X]$, if and only if
\[
\sum_{i=1}^n(u_i-v_i)\frac{\partial_Mf}{\partial_Mx_i}=0.
\]
The Jacobian matrix of an endomorphism $\varphi$ of $M(X)$ is
\[
J_M(\varphi)=\left(\frac{\partial_M\varphi(x_j)}{\partial_Mx_i}\right),
\]
which is a matrix with entries from $K[U,V]$.
One of the main results in \cite{U1} is that the
Jacobian matrix $J_M(\varphi)$ is invertible (as a matrix
with entries from $K[U,V]$) if and only if
$\varphi$ is an automorphism of $M(X)$.
Clearly, the invertibility of $J_M(\varphi)$ is equivalent
to $0\not= \text{\rm det}(J_M(\varphi))\in K$.
In this section we shall work with free algebras of rank 3 only and
shall assume that
the sets $X,U,V$ are, respectively,
\[
X=\{x,y,z\},\quad U=\{x_1,y_1,z_1\},\quad V=\{x_2,y_2,z_2\}.
\]
Let $T(K\langle x,y,z\rangle)$, $T(M(x,y,z))$ and $T(K[x,y,z])$ be,
respectively,
the groups of tame automorphisms of $K\langle x,y,z\rangle$,
$M(x,y,z)$, and
$K[x,y,z]$. There is a natural homomorphism
\[
T(K\langle x,y,z\rangle)\to T(M(x,y,z))\to T(K[x,y,z]).
\]
Let $\text{\rm Ker}(\pi)$ be the kernel of $\pi:T(M(x,y,z))\to
T(K[x,y,z])$.
Further developing the methodology  in \cite{SU1, SU2, SU3}, Umirbaev
\cite{U2}
discovered  the defining relations of $T(K[x,y,z])$. As a consequence
of that, he proved
the following.

\begin{proposition}
As a normal subgroup of $T(M(x,y,z))$, the kernel of $\pi$ is generated
by
the automorphisms
\[
\psi=(x+f(y,z),y,z),\quad f(y,z)
=\sum_{p,q,r,s\geq 0}\alpha_{pqrs}y^pz^q[y,z]y^rz^s,\,\alpha_{pqrs}\in
K.
\]
\end{proposition}

Moreover, any tame automorphism $\vartheta$ from $\text{\rm Ker}(\pi)$
has a Jacobian matrix
which is a product of elementary matrices.
The next key observation of Umirbaev is the following.
Let $\vartheta$ be any automorphism from the kernel of the
natural homomorphism $\text{Aut}M(x,y,z)\to \text{Aut}K[x,y,z]$.
Then
\[
J_M(\vartheta)=J_M(\vartheta)(x_1,y_1,z_1,x_2,y_2,z_2)
\]
is a $3\times 3$ matrix with entries from $K[x_1,y_1,z_1,x_2,y_2,z_2]$.
If we replace $x_1,y_1,x_2,y_2$ with zeros, then the matrix
$J_M(\vartheta)(0,0,z_1,0,0,z_2)$ will be of the form
\[
\overline{J_M(\vartheta)}(z_1,z_2)=\left(\begin{matrix}
1+w_{11}&w_{12}&w_{13}\\
w_{21}&1+w_{22}&w_{23}\\
0&0&1\\
\end{matrix}\right),
\]
where the polynomials $w_{ij}=w_{ij}(z_1,z_2)$ have no constant terms.
Define the $2\times 2$ matrix
\[
J_2(\vartheta)(z_1,z_2)=\left(\begin{matrix}
1+w_{11}(z_1,z_2)&w_{12}(z_1,z_2)\\
w_{21}(z_1,z_2)&1+w_{22}(z_1,z_2)\\
\end{matrix}\right).
\]

\begin{proposition}\label{proposition of Umirbaev} {\rm (Umirbaev
\cite{U2})}
If $\vartheta\in\text{\rm Ker}(\pi)$, then
$J_2(\vartheta)(z_1,z_2)$ is a product of elementary matrices with
entries
from $K[z_1,z_2]$.
\end{proposition}

Note that the matrix $J_2(\overline{\rho})$
of the automorphism $\overline{\rho}$ of $M(x,y,z)$ induced
by the automorphism $\rho$ of $K\langle x,y,z\rangle$ coincides with
the matrix
$J_z(\rho)$, the Jacobian matrix of $(\rho(x),\rho(y))$,
when $\rho$ fixes $z$ and is linear with respect to $x,y$.

Now we are ready to prove the main results in this article.

\begin{theorem}\label{general form of strong Anick conjecture}
Let $K$ be a field of characteristic $0$ and let
the polynomial $f(x,y,z)\in K\langle x,y,z\rangle$ be linear in $x,y$.
If there exists a wild automorphism of $K\langle x,y,z\rangle$ which
fixes $z$
and sends $x$ to $f(x,y,z)$, then every automorphism of $K\langle
x,y,z\rangle$ which
sends $x$ to $f(x,y,z)$ is also wild.
So, $f(x,y,z)$ is a wild coordinate of $K\langle x,y,z\rangle$.
\end{theorem}

\begin{proof}
Let $\sigma=(f(x,y,z),h(x,y,z),z)$ be a wild automorphism of $K\langle
x,y,z\rangle$
which fixes $z$ and sends $x$ to $f(x,y,z)$.
We write $f(x,y,z)$ in the form
\[
f(x,y,z)=\sum_{p,q\geq 0}\alpha_{pq}z^pxz^q+\sum_{p,q\geq
0}\beta_{pq}z^pyz^q+f_0(z),
\]
where $\alpha_{pq},\beta_{pq}\in K$,
and $f_0(z)$ is a polynomial in $z$. Let
\[
a(z_1,z_2)=\sum_{p,q\geq 0}\alpha_{pq}z_1^pz_2^q,
\]
\[
b(z_1,z_2)=\sum_{p,q\geq 0}\beta_{pq}z_1^pz_2^q.
\]
First we shall show that the polynomials $a(z_1,z_2),b(z_1,z_2)$ cannot
constitute the first
column of a matrix from $GE_2(K[z_1,z_2])$. Suppose on the contrary,
\[
J=\left(\begin{matrix}
a(z_1,z_2)&c(z_1,z_2)\\
b(z_1,z_2)&d(z_1,z_2)\\
\end{matrix}\right)\in GE_2(K[z_1,z_2]),
\]
\[
c(z_1,z_2)=\sum_{p,q\geq 0}\gamma_{pq}z_1^pz_2^q,
\]
\[
d(z_1,z_2)=\sum_{p,q\geq 0}\delta_{pq}z_1^pz_2^q.
\]
Consider the polynomial
\[
g(x,y,z)=\sum_{p,q\geq 0}\gamma_{pq}z^pxz^q+
\sum_{p,q\geq 0}\delta_{pq}z^pyz^q.
\]
By Proposition \ref{linear tame automorphisms fixing variable} (ii),
the automorphism $\rho=(f(x,y,z),g(x,y,z),z)$ is tame in the group of
automorphisms
fixing $z$.
Hence the automorphism
$\rho^{-1}\sigma$ is also wild. But
\[
\rho^{-1}\sigma=(x,k(x,y,z),z)
\]
for some $k(x,y,z)\in K\langle x,y,z\rangle$. This contradicts to Lemma
\ref{one nontrivial coordinate}.

Hence $a(z_1,z_2),b(z_1,z_2)$ cannot constitute the first
column of a matrix from $GE_2(K[z_1,z_2])$.

The next step is to produce a wild automorphism of $K\langle
x,y,z\rangle$ which induces the identity automorphism of
$K[x,y,z]$.

Let $\sigma=(f(x,y,z),h(x,y,z),z)$
be the above wild automorphism of $K\langle x,y,z\rangle$ which fixes
$z$
and sends $x$ to $f(x,y,z)$, and let $h_1(x,y,z)$ be the component
of $h$ which is linear with respect to $x,y$. Then
$\tau=(f(x,y,z),h_1(x,y,z),z)$ is also a wild automorphism of
$K\langle x,y,z\rangle$ which induces a $z$-tame automorphism of
$K[x,y,z]$. (The automorphism $\tau$ is wild since $\sigma^{-1}\tau$
sends $x$ to
$x$ and $z$ to $z$, then by Lemma \ref{one nontrivial coordinate}
$\sigma^{-1}\tau$ is tame.
The induced automorphism is $z$-tame by Proposition
\ref{linear tame automorphisms fixing variable}.)
Let $\psi$ be the corresponding $z$-tame, linear in $x,y$
automorphism of $K\langle x,y,z\rangle$. Then $\widetilde{\tau} =
\psi^{-1}\tau$ is still wild and induces the identity automorphism
of $K[x,y,z]$.

Now, let $\varphi$ be any tame automorphism of $K\langle
x,y,z\rangle$ which sends $x$ to $f(x,y,z)$. Replacing $\varphi$
with $\widetilde{\varphi}  = \psi^{-1}\varphi$, we obtain a tame
automorphism for which $\widetilde{\varphi}(x) =
\widetilde{\tau}(x)$.

The automorphism $\widetilde{\varphi}$ induces a tame automorphism
of $K[x,y,z]$ which fixes $x$. By  results in \cite{DY1, DY2, SU1,
SU3}, such an automorphism is tame in the class of automorphisms
fixing $x$ and we can lift it to an $x$-tame automorphism $\theta$
of $K\langle x,y,z\rangle$. So we obtain a tame automorphism
$\widehat{\varphi} = \widetilde{\varphi}\theta^{-1}$ which induces
the identity automorphism of $K[x,y,z]$ and $\widehat{\varphi}(x)
=  \widetilde{\tau}(x)$.

Let $\xi$ be the automorphism of $M(x,y,z)$ induced by
$\widehat{\varphi}$. It is in the kernel of the homomorphism $\pi$
of $\text{\rm Aut}M(x,y,z)\to \text{\rm Aut}K[x,y,z]$. The first
columns of the matrices $J_2(\xi)$ and
$J_2(\pi(\widetilde{\tau}))$ coincide. As we remarked above this
column cannot be a column of a matrix from $GE_2(K[z_1,z_2])$
since $\widetilde{\tau}$ is wild. On the other hand by
Proposition \ref{proposition of Umirbaev} it is a column of a
matrix from $GE_2(K[z_1,z_2])$. This contradiction completes the proof.
\end{proof}

Theorem \ref{general form of strong Anick conjecture} and
Proposition \ref{Euclidean algorithm for elementary matrices}
give an algorithm  deciding whether a polynomial
$f(x,y,z)\in K\langle x,y,z\rangle$
which is linear in $x$ and $y$, is a tame coordinate. If it is, then
the algorithm
finds a product of $z$-elementary automorphisms which sends $x$ to
$f(x,y,z)$.

The following consequence of Theorem
\ref{general form of strong Anick conjecture} gives the affirmative
answer to the
Strong Anick Conjecture.

\begin{theorem}
The Strong Anick Conjecture is true. Namely, there exist wild
coordinates
in $K\langle x,y,z\rangle$. In particular,
the two nontrivial coordinates $x+z(xz-zy)$ and $y+(xz-zy)z$ of the
Anick automorphism
\[
\omega=(x+z(xz-zy),y+(xz-zy)z,z)
\]
are both wild.
\end{theorem}

\begin{proof}
The partial derivatives of $f(x,y,z)=\omega(x)=x+z(xz-zy)$ are
\[
a(z_1,z_2)=\frac{\partial f}{\partial x}=1+z_1z_2,\quad
b(z_1,z_2)=\frac{\partial f}{\partial y}=-z_1^2.
\]
Since we cannot apply the Euclidean algorithm to $a(z_1,z_2)$ and
$b(z_1,z_2)$, Theorem
\ref{general form of strong Anick conjecture} gives that $f(x,y,z)$ is
a wild coordinate.
\end{proof}

We call an automorphism $\varphi=(f(x,y,z),g(x,y,z),z)$ of $K\langle
x,y,z\rangle$
{\it Anick-like} if $f(x,y,z)$ and $g(x,y,z)$ are linear in $x,y$ and
the matrix
$J_z(\varphi)$ does not belong to $GE_2(K[z_1,z_2])$.
The following corollary is an analogue of a result from \cite{UY}.

\begin{corollary}
The two nontrivial coordinates $f(x,y,z),g(x,y,z)$ of any Anick-like
automorphism
\[
\varphi=(f(x,y,z),g(x,y,z),z)
\]
of $K\langle x,y,z\rangle$ are wild.
\end{corollary}

\begin{proof}
Let
\[
\frac{\partial f}{\partial x}=a(z_1,z_2),\quad
\frac{\partial f}{\partial y}=b(z_1,z_2).
\]
We cannot apply the Euclidean algorithm to bring the pair
$(a(z_1,z_2),b(z_1,z_2))$
to $(\alpha,0)$, $0\not=\alpha\in K$,
because $J_z(\varphi)\not\in GE_2(K[z_1,z_2])$. Hence
Theorem
\ref{general form of strong Anick conjecture} gives that $f(x,y,z)$ is
a wild coordinate.
Similar arguments work for $g(x,y,z)$.
\end{proof}

In the spirit  of the above results, we obtain
the following theorem which is much stronger.

\begin{theorem}\label{more general theorem}
Let $f(x,y,z)$ be a $z$-coordinate of $K\langle x,y,z\rangle$
without terms depending only on $z$
(i.e. $f(0,0,z)=0$). If the linear part (with respect to $x$ and $y$)
$f_1(x,y,z)$ of $f(x,y,z)$
 is a $z$-wild coordinate,
then $f(x,y,z)$ itself is also a wild coordinate of $K\langle
x,y,z\rangle$.
\end{theorem}

\begin{proof} Since $f(x,y,z)$ is a $z$-coordinate of
$K\langle x,y,z\rangle$, there exists a $z$-automorphism
$\sigma=(f(x,y,z), g(x,y,z), z)$
of $K\langle x,y,z\rangle$. Obviously we may assume
$g(0,0,z)=0$ (otherwise just replace $g(x,y,z)$
by ($g(x,y,z)-g(0,0,z)$). Let $\sigma_1=(f_1(x,y,z),g_1(x,y,z),z)$ be
the automorphism
which is the linear part of $\sigma$. By assumption $\sigma_1$ is a
wild automorphism.
We have to prove the wildness of all automorphisms
$\varphi=(f(x,y,z),u(x,y,z),v(x,y,z))$ of $K\langle x,y,z\rangle$ with
first coordinate equal to $f(x,y,z)$.
Consider the automorphisms
$\overline{\sigma}=(\overline{f},\overline{g},z)$
and $\overline{\varphi}=(\overline{f},\overline{u},\overline{v})$ of
$K[x,y,z]$
induced by $\sigma$ and $\varphi$, respectively. If $\overline{\sigma}$
is wild, then,
by the theorem of Umirbaev and Yu \cite{UY}, $\overline{f}$ is a wild
coordinate of
$K[x,y,z]$. Hence $\overline{\varphi}$ is a wild automorphism of
$K[x,y,z]$. This implies that
$\varphi$ is a wild automorphism of $K\langle x,y,z\rangle$ and
therefore
$f(x,y,z)$ is a wild coordinate.
Hence we may assume that $\overline{\sigma}$ is a tame automorphism of
$K[x,y,z]$. Now we suppose that the automorphism $\varphi$ is tame and
repeat the main steps of the proof of Theorem
\ref{general form of strong Anick conjecture}. Since
$\overline{\sigma}$ is tame,
by \cite{DY1, DY3, SU1, SU3} it is also $z$-tame. Let $\psi$ be some
$z$-tame automorphism of $K\langle x,y,z\rangle$ which induces
$\overline{\sigma}$
and let $\psi_1$ be the linear part of $\psi$.
Replacing $\sigma$ with $\widetilde{\sigma}=\psi^{-1}\sigma$ and
$\varphi$ with $\widetilde{\varphi}=\psi^{-1}\varphi$, we obtain that
the tame automorphism
$\widetilde{\varphi}$ fixes $x$ modulo the commutator ideal of
$K\langle x,y,z\rangle$.
Since $\widetilde{\sigma}$ is a composition of the $z$-automorphisms
$\psi^{-1}$ and $\sigma$,
its linear part $(\widetilde{\sigma})_1$ is also a $z$-automorphism
which is equal to the composition
$\psi_1^{-1}\varphi_1$ of the linear components of $\psi_1^{-1}$ and
$\varphi_1$.
Hence $(\widetilde{\sigma})_1$ is wild and we may reduce our
considerations to the case when
$\widetilde{\sigma}(x)=f(x,y,z)$
is congruent to $x$ modulo the commutator ideal of $K\langle
x,y,z\rangle$.
Since $\widetilde{\varphi}$ induces a tame automorphism of $K[x,y,z]$,
by \cite{DY1, DY3, SU1, SU3}
again, the induced automorphism is also $x$-tame and we can lift it to
an $x$-tame automorphism
$\theta$ of $K\langle x,y,z\rangle$. The tame automorphism
$\widehat{\varphi} = \widetilde{\varphi}\theta^{-1}$ induces
the identity automorphism of $K[x,y,z]$ and
$\widehat{\varphi}(x)=f(x,y,z)$.
Now, as in Theorem \ref{general form of strong Anick conjecture},
the proof is completed with considerations in the free metabelian
algebra $M(x,y,z)$.
\end{proof}

\begin{remark}
The restriction $f(0,0,z)=0$ is essential for the proof of
Theorem \ref{more general theorem} (Note that obviously we may assume
$g(0,0,z)=0$, otherwise just replace $g(x,y,z)$
by $g(x,y,z)-g(0,0,z)$). We use it when,
modifying simultaneously the automorphisms
$\sigma=(f(x,y,z), g(x,y,z), z)$ and
$\varphi=(f(x,y,z),u(x,y,z),v(x,y,z))$
of $K\langle x,y,z\rangle$, we bring $\sigma$ and $\varphi$ to
automorphisms which send $x$
to the same element congruent
to $x$ modulo the commutator ideal, still keeping the property that the
linear component
of the image of $x$ is wild. Nevertheless, it seems very unlikely to
have a wild automorphism
$(f,g,z)$ with $f(0,0,z)=0$ such that $f(x,y,z)+a(z)$ is a
tame coordinate
for some polynomial $a(z)$ in view of the next theorem.
\end{remark}

\begin{theorem}\label{nonlinear automorphisms with z fixed}
Let $(f,g,z)$ be an automorphism of $K\langle x,y,z\rangle$
and let the linear
part (with respect to $x$ and $y$) of it, $(f_1,g_1,z)$,  be a $z$-wild
automorphism.
Then $(f,g,z)$ is also a wild automorphism of $K\langle x,y,z\rangle$.
\end{theorem}

\begin{proof}
Let $f(x,y,z)=f'(x,y,z)+f_0(z)$, $g(x,y,z)=g'(x,y,z)+g_0(z)$,
where $f',g'$ do not contain monomials depending on $z$ only.
Define the automorphism $\tau=(x-f_0(z),y-g_0(z),z)$. Then the
automorphism
$\sigma=(f,g,z)$ is tame (or $z$-tame) if and only if
$\sigma\tau=(f',g',z)$ is tame (or $z$-tame). Since the polynomials
$f,f'$ and
$g,g'$ have the same linear components $f_1$ and $g_1$, we apply
Theorem \ref{more general theorem}.
\end{proof}

\begin{remark} The above theorem is much stronger
than the main result in \cite{U2} where only the
automorphisms linear with respect to $x$  and $y$ are dealt.
\end{remark}

The following example gives a large class of wild automorphisms
and wild coordinates. It is based on the polynomial $xz-zy$ which
appears in the Anick automorphism.

\begin{example}\label{example generalizing Anick}
Let $h(t,z)\in K\langle t,z\rangle$ and let $h(0,0)=0$.
Then
\[
\sigma_h=(x+zh(xz-zy,z),y+h(xz-zy,z)z,z)
\]
is an automorphism of $K\langle x,y,z\rangle$  fixing $xz-zy$.
If the linear component (with respect to $x,y$)
$h_1(xz-zy,z)$ of $h(xz-zy,z)$ is not equal to 0,
then this automorphism belongs to the class of wild automorphisms
in Theorem \ref{more general theorem}: As
$(\sigma_h)_1=(x+zh_1(xz-zy,z),y+h_1(xz-zy,z)z,z)$
is an automorphism of $K\langle x,y,z\rangle$
and its matrix $J_z((\sigma_h)_1)$ is
\[
J_z((\sigma_h)_1)=\left(\begin{matrix}
1+q(z_1,z_2)z_1z_2&q(z_1,z_2)z_2^2\\
-q(z_1,z_2)z_1^2&1-q(z_1,z_2)z_1z_2\\
\end{matrix}\right)
\]
for some nonzero polynomial $q(z_1,z_2)\in K[z_1,z_2]$,
it is easy to see that this matrix does not belong to
$GL_2(K[z_1,z_2])$
because we cannot apply the Euclidean algorithm to its first column.
\end{example}

\begin{example}\label{example without linear terms}
A minor modification of the Anick automorphism is the automorphism
of $K\langle x,y,z\rangle$
\[
\omega_m=(x+z(xz-zy)^m,y+(xz-zy)^mz,z).
\]
Note that the automorphisms $\omega_m$, $m>1$,
are not covered by Theorem \ref{more general theorem}, as
the polynomials $z(xz-zy)^m$ and $(xz-zy)^mz$ have no
linear components with respect to $x$ and $y$.
\end{example}

\begin{theorem}\label{generalization of Anick automorphism}
The above automorphisms $\omega_m$ are wild for all $m\geq 1$.
\end{theorem}

\begin{proof}
Consider the automorphism $\tau=(x+1,y,z)$ of $K\langle x,y,z\rangle$.
Clearly, $\omega_m$ is wild if and only if $\omega_m\tau$ is wild.
Direct calculations show that the linear part of the $z$-automorphism
\[
\omega_m\tau=x+1+z((x+1)z-zy)^m,y+((x+1)z-zy)^mz,z)
\]
is equal to
\[
(\omega_m\tau)_1=\left(x+z\sum_{i=0}^{m-1}z^i(xz-zy)z^{m-1-i},
y+\sum_{i=0}^{m-1}z^i(xz-zy)z^{m-1-i}z,z\right).
\]
Hence the matrix $J_z((\omega_m\tau)_1)$ has the form
\[
J_z((\omega_m\tau)_1)=\left(\begin{matrix}
1+q(z_1,z_2)z_1z_2&q(z_1,z_2)z_2^2\\
-q(z_1,z_2)z_1^2&1-q(z_1,z_2)z_1z_2\\
\end{matrix}\right),
\]
where $q(z_1,z_2)=z_1^{m-1}+z_1^{m-2}z_2+\cdots+z_2^{m-1}$.
As in Example \ref{example generalizing Anick},
the automorphism $(\omega_m\tau)_1$ is wild. Hence $\omega_m$ is also
wild
by Corollary \ref{nonlinear automorphisms with z fixed}.
\end{proof}

It seems plausible that the nontrivial coordinates of $\omega_m$,
$m>1$, are wild.
However, our methods and the the methods in \cite{U2}
are not applicable here.

\begin{problem}
Are the two nontrivial coordinates of the above automorphism
$\omega_m$, $m>1$, both wild?
\end{problem}

\begin{remark}
The most general form of the result of
Umirbaev \cite{U2} gives that the automorphism $\vartheta=(f,g,h)$ of
the free metabelian algebra $M(x,y,z)$ is wild, if
it induces the identity automorphism of $K[x,y,z]$
and the matrix $J_2(\vartheta)(z_1,z_2)$
cannot be presented as a product of elementary matrices with entries
from $K[z_1,z_2]$, see Proposition \ref{proposition of Umirbaev}.
Hence the classes of wild automorphisms and wild coordinates in
Theorem \ref{more general theorem}, Example \ref{example generalizing
Anick}
and Example \ref{example without linear  terms} are
not covered by Umirbaev \cite{U2}.
\end{remark}

Now we are going to show that at least two coordinates of the
automorphisms
of the class of Umirbaev are wild.

\begin{theorem}\label{the case of the theorem of Umirbaev}
Let $\vartheta=(f,g,h)$ be an automorphism of
the free metabelian algebra $M(x,y,z)$ which
induces the identity automorphism of $K[x,y,z]$
and the matrix $J_2(\vartheta)(z_1,z_2)$ does not belong to
$GE_2(K[z_1,z_2])$.
Then the two  coordinates $f(x,y,z)$ and $g(x,y,z)$ are both wild.
\end{theorem}

\begin{proof}
We repeat the main steps of the proof of
Theorem \ref{general form of strong Anick conjecture}.
The polynomial $f(x,y,z)\in M(x,y,z)$ is equal to $x$ modulo
the commutator ideal of $M(x,y,z)$ and has the form
\[
f=\sum_{p,q\geq 0}\alpha_{pq}z^pxz^q+\sum_{p,q\geq 0}\beta_{pq}z^pyz^q+
\sum_{k\geq 2}f_k(x,y,z),
\]
where $f_i$ is the homogeneous component of degree $i$ in $x,y$ (and
$f_0=0$). Let
\[
a(z_1,z_2)=\sum_{p,q\geq 0}\alpha_{pq}z_1^pz_2^q,\quad
b(z_1,z_2)=\sum_{p,q\geq 0}\beta_{pq}z_1^pz_2^q.
\]
The polynomials $a(z_1,z_2),b(z_1,z_2)$ constitute the first column of
the matrix
$J_2(\vartheta)(z_1,z_2)$ which does not belong to
$GE_2(K[z_1,z_2])$. By Proposition \ref{Euclidean algorithm for
elementary matrices}
$a(z_1,z_2),b(z_1,z_2)$ cannot be reduced to $(\alpha,0)$,
$0\not=\alpha\in K$, by the
Euclidean algorithm only.

Now, let $\varphi=(f(x,y,z),u(x,y,z),v(x,y,z))$ be any tame
automorphism which sends
$x$ to $f(x,y,z)$. Clearly, $\varphi$ induces the tame automorphism
\[
\overline{\varphi}=(\overline{f},\overline{u},\overline{v})=(x,\overline{u},\overline{v})
\]
of the polynomial algebra $K[x,y,z]$. Since $\overline{\varphi}$ fixes
$x$,
the results in  \cite{DY1, DY2, SU1, SU2, SU3} give that
$\overline{\varphi}$ is tame also in the class of automorphisms fixing
$x$.
So, as in the proof of Theorem \ref{general form of strong Anick
conjecture},
we may replace $\varphi$ with a tame automorphism
$\xi=(f(x,y,z),u_1(x,y,z),v_1(x,y,z))$ of $M(x,y,z)$ such that
$\xi$ is in the kernel of the natural homomorphism
$\text{Aut}M(x,y,z)\to \text{Aut}K[x,y,z]$.
The tameness of $\xi$ implies that $J_2(\xi)\in GE_2(K[z_1,z_2])$.
Since the first column of
$J_2(\xi)$ consists of $a(z_1,z_2)$ and $b(z_1,z_2)$, this contradicts
to
Proposition \ref{proposition of Umirbaev}. The considerations for the
other coordinate
$g$ of $\vartheta$ are similar.
\end{proof}

\begin{remark}
Any automorphism $\phi\in\text{Aut}K\langle x,y,z\rangle$
which induces an automorphism in $\text{Aut}M(x,y,z)$ of the type in
Theorem \ref{the case of the theorem of Umirbaev}
(in other words, any automorphism in $\text{Aut}K\langle x,y,z\rangle$
obtained by lifting an automorphism in $\text{Aut}M(x,y,z)$ of the type
in
Theorem \ref{the case of the theorem of Umirbaev})
is a wild automorphism  containing at least two wild coordinates.
\end{remark}

The above results suggest the following problems.

\begin{problem}
Is it true that the two nontrivial coordinates of
a wild automorphism of $K\langle x,y,z\rangle$  fixing $z$ are both
wild?
\end{problem}

\begin{problem}
Is it true that every wild automorphism of $K\langle x,y,z\rangle$
contains at least two wild coordinates?
\end{problem}

\section{Special wild automorphisms of the free metabelian algebra}

In this section we shall construct a wild automorphism $\tau$ of the
free metabelian
algebra $M(x,y,z)$ over any field $K$ of arbitrary characteristic
with the following properties:

(i) $\tau=(f(x,y,z),y,z)$
fixes two of the variables. (Hence Lemma \ref{one nontrivial
coordinate} does not hold
for $M(x,y,z)$.)

(ii) The Jacobian matrix $J_M(\tau)$ is a product of elementary
matrices.

(iii) It cannot be lifted to an automorphism of $K\langle
x,y,z\rangle$.

Recall the definition of the Fox derivatives of the free algebra
$K\langle X\rangle$,
see e.g. \cite{MSY}.
If
\[
f(X)=\sum_{i=1}^nx_if_i(X)+\alpha,\quad
\alpha\in K,\quad f_i(X)\in K\langle X\rangle,
\]
then the {\it right Fox derivatives} of $f(X)$ are
\[
\frac{\partial_r f}{\partial_rx_i}=f_i(X),\quad i=1,\ldots,n.
\]
Similarly, if
\[
f(X)=\sum_{i=1}^nf_i(X)x_i+\alpha,\quad
\alpha\in K,\quad f_i(X)\in K\langle X\rangle,
\]
then the {\it left Fox derivatives} of $f(X)$ are
\[
\frac{\partial_l f}{\partial_lx_i}=f_i(X),\quad i=1,\ldots,n.
\]
The right and left Jacobian matrices of an endomorphism $\varphi$ of
$K\langle X\rangle$ are, respectively,
\[
J_r(\varphi)=\left(\frac{\partial_r\varphi(x_j)}{\partial_rx_i}\right),\quad
J_l(\varphi)=\left(\frac{\partial_l\varphi(x_j)}{\partial_lx_i}\right).
\]
The chain rule gives that if $\varphi$ is an automorphism, then
$J_r(\varphi)$ and $J_l(\varphi)$ are invertible (but the oposite is
not true in the
general case).

We need some machinery from \cite{BGLM} and \cite{BD}.
We state it in the case of three variables only. We define an
equivalence
relation $\sim$ on $K\langle x,y,z\rangle$. We say that two monomials
$u$ and $v$
are equivalent, if they can be obtained from each other by cyclic
permutation
(i.e., $u\sim v$ if and only if $u=w_1w_2$ and $v=w_2w_1$ for some
monomials $w_1,w_2$)
and then extend $\sim$ to $K\langle x,y,z\rangle$ by linearity.

\begin{proposition}\label{trace test}
\cite{BGLM, BD}
Let $\sigma$ be an endomorphism of $K\langle x,y,z\rangle$
which is equal to the identity of $K\langle x,y,z\rangle$
modulo the $k$-th degree of the augmentation ideal, i.e.
\[
\sigma=(x+f_k+\cdots+f_m,y+g_k+\cdots+g_m,z+h_k+\cdots+h_m),
\]
where $f_i,g_i,h_i$ are the homogeneous
components of degree $i$ of $\sigma(x),\sigma(y),\sigma(z)$,
respectively.
If $\sigma$ is an automorphism and $k\geq 2$, then the homogeneous
component of
degree $k-1$ of the trace of the right Jacobian matrix
\[
\frac{\partial_rf_k}{\partial_rx}
+\frac{\partial_rg_k}{\partial_ry}+\frac{\partial_rh_k}{\partial_rz}
\]
is equivalent to $0$. Similar statement holds for the trace of the left
Jacobian matrix.
\end{proposition}

\begin{theorem}
The endomorphism
\[
\tau=(x+x^2[y,z],y,z)
\]
of the free metabelian algebra $M(x,y,z)$ is a wild automorphism which
cannot be lifted to
an automorphism of $K\langle x,y,z\rangle$. Its Jacobian matrix
\[
J_M(\tau)=\left(\begin{matrix}
1&0&0\\
x_1^2(z_2-z_1)&1&0\\
x_1^2(y_1-y_2)&0&1\\
\end{matrix}\right)
\]
is a product of two elementary matrices.
\end{theorem}

\begin{proof}
Obviously $\tau$ is an automorphism and $\tau^{-1}=(x-x^2[y,z],y,z)$.
Also, its Jacobian matrix $J_M(\tau)$ is a product of elementary
matrices.
Now, let $\tau$ be lifted to an automorphism $\sigma$ of $K\langle
x,y,z\rangle$.
Then
\[
\sigma=(x+x^2[y,z]+f(x,y,z),y+g(x,y,z),z+h(x,y,z)),
\]
where $f(x,y,z),g(x,y,z),h(x,y,z)$ belong to the T-ideal generated by
the polynomial identity $[x_1,x_2][x_3,x_4]=0$. Hence
$f,g,h$ have no homogeneous components of degree $\leq 3$ and
\[
f=f_4+\cdots+f_m,\quad g=g_4+\cdots+g_m,\quad h=h_4+\cdots+h_m,
\]
where $f_i,g_i,h_i$ are homogeneous of degree $i$. Clearly,
the components $f_4,g_4,h_4$ are linear combinations of products
of two commutators of the variables. By Proposition \ref{trace test},
\begin{equation}\label{test for right derivatives}
\frac{\partial_r(x^2[y,z]+f_4)}{\partial_rx}+
\frac{\partial_rg_4}{\partial_ry}+\frac{\partial_rh_4}{\partial_rz}\sim
0,
\end{equation}
\begin{equation}\label{test for left derivatives}
\frac{\partial_l(x^2[y,z]+f_4)}{\partial_lx}+
\frac{\partial_lg_4}{\partial_ly}+\frac{\partial_lh_4}{\partial_lz}\sim
0.
\end{equation}
Since $x^2[y,z]=x^2yz-x^2zy$, we obtain that
\begin{equation}\label{Fox derivatives for sigma}
\frac{\partial_rx^2[y,z]}{\partial_rx}=x[y,z]\sim xyz-xzy, \quad
\frac{\partial_lx^2[y,z]}{\partial_lx}=0.
\end{equation}
The components of
(\ref{test for right derivatives}) and (\ref{test for left
derivatives})
which are multilinear in $x,y,z$ are equivalent to 0.
The components of the Fox derivatives
\[
\frac{\partial_rf_4}{\partial_rx},\quad
\frac{\partial_rg_4}{\partial_ry},\quad
\frac{\partial_rh_4}{\partial_rz}
\]
which are multilinear in $x,y,z$ come, respectively, from
\[
f_4'=\alpha_1[x,y][x,z]+\beta_1[x,z][x,y],
\]
\[
g_4'=\alpha_2[x,y][y,z]+\beta_2[y,z][x,y],
\]
\[
h_4'=\alpha_3[x,z][y,z]+\beta_3[y,z][x,z].
\]
Direct calculations give that
\[
\frac{\partial_rf_4'}{\partial_rx}+
\frac{\partial_rg_4'}{\partial_ry}+\frac{\partial_rh_4'}{\partial_rz}
\]
\[
\sim (\alpha_1y[x,z]+\beta_1z[x,y])+(-\alpha_2x[y,z]+\beta_2z[x,y])
-(\alpha_3x[y,z]+\beta_3y[x,z])
\]
\[
\sim (-\alpha_1+\beta_1-\alpha_2+\beta_2-\alpha_3+\beta_3)(xyz-xzy).
\]
Together with (\ref{Fox derivatives for sigma}) this implies that
\begin{equation}\label{right trace equal to zero}
-\alpha_1+\beta_1-\alpha_2+\beta_2-\alpha_3+\beta_3+1=0.
\end{equation}
Similarly,
\[
\frac{\partial_lf_4'}{\partial_lx}+
\frac{\partial_lg_4'}{\partial_ly}+\frac{\partial_lh_4'}{\partial_lz}
\]
\[
\sim -(\alpha_1[x,y]z+\beta_1[x,z]y)+(-\alpha_2[x,y]z+\beta_2[y,z]x)
+(\alpha_3[x,z]y+\beta_3[y,z]x)
\]
\[
\sim (-\alpha_1+\beta_1-\alpha_2+\beta_2-\alpha_3+\beta_3)(xyz-xzy)\sim
0
\]
in virtue of (\ref{Fox derivatives for sigma}). Hence
\begin{equation}\label{left trace equal to zero}
-\alpha_1+\beta_1-\alpha_2+\beta_2-\alpha_3+\beta_3=0.
\end{equation}
Clearly, (\ref{right trace equal to zero}) and (\ref{left trace equal
to zero})
contradict to each other. Hence $\tau$ cannot be lifted to an
automorphism
of $K\langle x,y,z\rangle$ and, therefore, is a wild automorphism of
$M(x,y,z)$.
\end{proof}

\begin{problem}
Is the polynomial $x+x^2[y,z]$ a wild coordinate of $M(x,y,z)$?
Can it be lifted to a coordinate of $K\langle x,y,z\rangle$?
\end{problem}

\begin{problem}
Do there exist wild automorphisms and wild coordinates
of the free metabelian algebra $M(X)$
of rank $n>3$? Are there wild automorphisms similar to the automorphism
$\tau$ constructed above?
\end{problem}

\section{Lifting of automorphisms fixing variables}

The considerations in this section work over an arbitrary field of any
characteristic.

Let $G(X)$ be the free group generated by the finite set $X$.
The theorem of Nielsen \cite{Ni} states that every automorphism of
$G(X)$
is a product of the elementary automorphisms
$(x_1^{-1},x_2,\ldots,x_n)$, $(x_1x_2,x_2,\ldots,x_n)$, and
$(x_{\sigma(1)},\ldots,x_{\sigma(n)})$,
where $\sigma$ belongs to the symmetric group $S_n$.
The proof of Nielsen gives
also an algorithm which finds such a decomposition. The theorem of
Schreier
\cite{Sch} states that every subgroup of the free group with any number
of generators
is also free.

There are several important varieties of algebras over a field with
free objects
which share the above properties of free groups. The variety $\mathfrak
V$
of algebras is called {\it Schreier} if the subalgebras of the
relatively free
algebras $F({\mathfrak V})$ are again relatively free,
where $F({\mathfrak V})$ is freely generated by a set of any
cardinality.
The variety $\mathfrak V$ is {\it Nielsen} if all automorphisms of the
free algebras
$F_n({\mathfrak V})$ of finite rank are tame. A theorem of Lewin
\cite{L}
gives that over an infinite field $K$ the two notions coincide, i.e.,
$\mathfrak V$ is Nielsen if and only if it is Schreier. The same holds
over an
arbitrary field $K$, provided that the variety
$\mathfrak V$ is defined by a multilinear system of polynomial
identities.
See the book \cite{MSY} for more details about examples of Schreier
varieties,
and the properties of the subalgebras and the automorphisms of their
free objects.

The variety of all (not necessarily associative)
algebras is Schreier, by the theorem of Kurosh \cite{Ku}.
Recall that the absolutely free algebra $K\{X\}$ consists of all
polynomials in the set of noncommuting and nonassociative variables
$X$, e.g.
$(xx)x\not=x(xx)$.
One of the key moments of the proof of Kurosh (and of all other proofs
that some varieties are Schreier) is the following, see \cite{MSY},
Theorem 11.1.1.
For a nonzero polynomial $f\in K\{X\}$ we denote by $\bar f$
the homogeneous component of maximal degree of $f$.

\begin{proposition}\label{theorem of Kurosh}
{\rm (i)} Any finite set $S$ of $K\{X\}$ can be transformed into a set
of free
generators of the subalgebra generated by $S$ by a finite sequence of
elementary
transformations (with cancellation of possible zeros).

{\rm (ii)} If $F=\{f_1,\ldots,f_n\}$ is a set of free generators of
$K\{X\}$,
and $g\in K\{X\}$, then $\bar g$ belongs to the subalgebra of $K\{X\}$
generated by $\overline{f_1},\ldots,\overline{f_n}$.
\end{proposition}

For an automorphism $\varphi=(f_1,\ldots,f_n)$ of $K\{X\}$ we define
the degree of $\varphi$ as the sum of the degrees of the coordinates
$f_i$:
\[
\text{deg}(\varphi)=\sum_{i=1}^n\text{deg}(f_i).
\]
Clearly, $\text{deg}(\varphi)\geq n$. The following consequence of
Proposition \ref{theorem of Kurosh} can be used
effectively to decompose an automorphism of $K\{X\}$
as a product of elementary automorphisms.

\begin{corollary}\label{decoposition of automorphisms of free algebras}
Let $\varphi=(f_1,\ldots,f_n)\in\text{\rm Aut}K\{X\}$ with
$\text{\rm deg}(\varphi)>n$. Then there exists an integer $i$ and a
polynomial
$g(y_1,\ldots,y_{i-1},y_{i+1},\ldots,y_n)$, such that
\[
\overline{f_i}=g(\overline{f_1},\ldots,\overline{f_{i-1}},
\overline{f_{i+1}},\ldots,\overline{f_n}).
\]
Let $\tau$ be the elementary automorphism of $K\{X\}$ defined by
\[
\tau=(x_1,\ldots,x_{i-1},x_i-g(x_1,\ldots,x_{i-1},x_{i+1},\ldots,x_n),
x_{i+1},\ldots,x_n).
\]
Then $\text{\rm deg}(\varphi\tau)<\text{\rm deg}(\varphi)$.
\end{corollary}

Now we are able to prove the following.

\begin{theorem}\label{Z-automorphisms of free algebras}
Let
\[
\varphi=(f_1(X,Z),\ldots,f_n(X,Z),z_1,\ldots,z_m)\in\text{\rm
Aut}_ZK\{X,Z\}
\]
be an automorphism of $K\{X,Z\}$ fixing the variables $Z$. Then
$\varphi$
is tame in the class of $Z$-automorphisms.
\end{theorem}

\begin{proof}
Let us consider $\varphi$ as an automorphism of $K\{X,Z\}$ in the usual
sense. The total degree of $\varphi$ is
\[
\text{\rm deg}(\varphi)=\sum_{i=1}^n\text{\rm deg}(f_i(X,Z))
+\sum_{j=1}^m\text{\rm deg}(z_j)
=\sum_{i=1}^n\text{\rm deg}(f_i)+m.
\]
Since $\varphi$ is a $Z$-automorphism, each polynomial $f_1,\ldots,f_n$
essentially depends on $X$.

If $\text{\rm deg}(\varphi)=n+m$, then all polynomails $f_i(X,Z)$ are
of total degree
equal to 1
and $\varphi$ is affine. We replace $\varphi$ with the product
$\psi=\varphi\tau_0$, where $\tau_0$ is the translation
\[
\tau_0=(x_1-f_1(0,0),\ldots,x_n-f_n(0,0),z_1,\ldots,z_m).
\]
Clearly, $\tau_0$ is a product of $Z$-elementary automorphisms
and $\psi$ is a linear automorphism. Its matrix, as a linear operator
of the vector space with basis $X\cup Z$, is
\[
\left(\begin{matrix}
A&0\\
B&E_m\\
\end{matrix}\right)
=\left(\begin{matrix}
\alpha_{11}&\alpha_{12}&\ldots&\alpha_{1n}&0&0\ldots&0\\
\alpha_{21}&\alpha_{22}&\ldots&\alpha_{2n}&0&0\ldots&0\\
\vdots&\vdots&\ddots&\vdots&\vdots&\vdots\ddots&\vdots\\
\alpha_{n1}&\alpha_{n2}&\ldots&\alpha_{nn}&0&0\ldots&0\\
\beta_{11}&\beta_{12}&\ldots&\beta_{1n}&1&0\ldots&0\\
\beta_{21}&\beta_{22}&\ldots&\beta_{2n}&0&1\ldots&0\\
\vdots&\vdots&\ddots&\vdots&\vdots&\vdots\ddots&\vdots\\
\beta_{m1}&\beta_{m2}&\ldots&\beta_{mn}&0&0\ldots&1\\
\end{matrix}\right),
\]
and $A=(\alpha_{pq})$, $B=(\beta_{rs})$ are, respectively,
$n\times n$ and $m\times n$ matrices with entries in
$K$, $E_m$ is the $m\times m$ identity matrix, and $A$ is invertible.
Since we work over a field, $A$ is a product of elementary matrices
and this implies that, multiplying $\psi$ by a product of elementary
linear
automorphisms fixing $Z$, we bring it to the automorphism
\[
\tau_1=(x_1+g_1(Z),\ldots,x_n+g_n(Z),z_1,\ldots,z_m),
\]
which is a product of elementary automorphisms fixing $Z$.

Now, let $\text{\rm deg}(\varphi)>n+m$. Then, at least one of the
polynomials
$f_i(X,Z)$ is not linear. The leading components of the $n+m$
coordinates are
$$\overline{f_1(X,Z)},\ldots, \overline{f_n(X,Z)},\
\overline{z_1}=z_1,\ldots,\overline{z_m}=z_m.$$
By Corollary \ref{decoposition of automorphisms of free algebras},
one of these homogeneous components is expressed by the others.
Obviously,
$z_j$ cannot be expressed as a polynomial of
$\overline{f_1},\ldots,\overline{f_n}$
and the other $z_1,\ldots,z_{j-1},z_{j+1},\ldots,z_m$. Hence, some
$\overline{f_i}$
is a polynomial of
$\overline{f_1},\ldots,\overline{f_{i-1}},\overline{f_{i+1}},
\ldots,\overline{f_n}$ and $z_1,\ldots,z_m$. This gives that the
elementary automorphism $\tau$ of $K\{X,Z\}$ prescribed by
Corollary \ref{decoposition of automorphisms of free algebras}, is of
the form
\[
\tau=(x_1,\ldots,x_{i-1},x_i-g(X,Z),
x_{i+1},\ldots,x_n,z_1,\ldots,z_m),
\]
where $g(X,Z)$ does not depend on $x_i$.
Then $\text{\rm deg}(\varphi\tau)<\text{\rm deg}(\varphi)$ and the
proof is completed by
obvious induction on the degree of $\varphi$.
\end{proof}

The theorem below is an immediate concequence of Theorem
\ref{Z-automorphisms of free algebras}.

\begin{theorem}\label{wild Z-automorphisms cannot be lifted to
Z-automorphisms}
Let $\varphi$ be an automorphism of $K[X,Z]$ or $K\langle X,Z\rangle$
which fixes $Z$.
If $\varphi$ is wild as a $Z$-automorphism, then it cannot be lifted to
an automorphism
of $K\{X,Z\}$ which also fixes $Z$.
\end{theorem}

\begin{remark} The Nagata automorphism is wild as a $z$-automorphism of
$K[x,y,z]$, \cite{N}, as well as wild in the usual sense, \cite{SU1,
SU3}.
Hence it cannot be lifted to any automorphism of $K\{x,y,z\}$.
On the other hand, by the theorem of Smith \cite{Sm},
automorphisms of $K[X]$ of a large class become tame as automorphisms
of $K[X,t]$,
if we extend them to act identically $t$. In particular, the extension
of
the Nagata automorphism
\[
\nu'=(x-2y(y^2+xz)-z(y^2+xz)^2,y+z(y^2+xz),z,t)
\]
is tame as an automorphism of $K[x,y,z,t]$. It is easy to see that it
is
wild in the class of automorphisms of $K[x,y,z,t]$ fixing $z$ and $t$.
Hence, Theorem
\ref{wild Z-automorphisms cannot be lifted to Z-automorphisms} gives
that
$\nu'$ cannot be lifted to an automorphism of $K\{x,y,z,t\}$ which
fixes $z$ and $t$.

Similarly, the automorphism of Anick is wild as an automorphism fixing
a variable
\cite{DY3} and even wild in the usual sense \cite{U2}. But it becomes
tame
extended to an automorphism of $K\langle x,y,z,t\rangle$. The technique
of
\cite{DY3} gives that the extension of the Anick automorphism
\[
(x+z(xz-zy),y+(xz-zy)z,z,t)
\]
is wild in the group of automorphisms of $K\langle x,y,z,t\rangle$
which fix
$z,t$, although this automorphism is tame in the usual sense.
Hence, our theorem gives that it cannot be lifted to an automorphism
of $K\{x,y,z,t\}$ which fixes $z,t$.
\end{remark}

We shall conclude this section with several open problems.

\begin{problem}
{\rm (i)} If $\varphi$ is an automorphism of $K[X]$, can it be
lifted to an automorphism of $K\langle X\rangle$? (If $\varphi$ is
wild, and nevertheless the answer is positive, this would mean
that it is not ``too wild''.)

{\rm (ii)} If $\varphi\in \text{\rm Aut}_ZK[X,Z]$, can it be
lifted to a $Z$-automorphism of $K\langle X,Z\rangle$? Are
$Z$-wild automorphisms wild also in the usual sense?
\end{problem}

\begin{problem} How far can be lifted the automorphisms of $K[X]$?
Describe the varieties of algebras $\mathfrak V$ with the property
that every automorphism of $K[X]$ can be lifted to an automorphism
of the relatively free algebra $F_n({\mathfrak V})$ of rank
$n=\vert X\vert$.
\end{problem}

For example, a theorem of Umirbaev \cite{U1} gives that every
automorphism of $K[X]$ can be lifted to an automorphism of the
free metabelian algebra $M(X)$.

\begin{problem}
{\rm (i)} If $p(X)$ is a coordinate of $K[X]$, can it be lifted
to a coordinate of $K\langle X\rangle$?

{\rm (ii)} How far can be lifted the coordinates of $K[X]$?
Describe the varieties of algebras $\mathfrak V$ with the property
that every coordinate of $K[X]$ can be lifted to a coordinate of
$F_n({\mathfrak V})$.

{\rm (iii)} Can the two nontrivial Nagata coordinates
$x-2y(y^2+xz)-z(y^2+xz)^2$ and $y+z(y^2+xz)$ be lifted to coordinates
of $K\langle x,y,z\rangle$?
\end{problem}

\section*{Acknowledgements}

The authors are grateful to Leonid Makar-Limanov, Ivan Shestakov
and Efim Zelmanov for insightful discussion.

\end{document}